\documentclass{article}%
\usepackage{amsmath}%
\usepackage{amsfonts}%
\usepackage{amssymb}%
\usepackage{graphicx}

\begin{document}

\title{On Transforming the Generalized Exponential Power Series}
\author{Henrik Stenlund\thanks{The author is grateful to Visilab Signal Technologies for supporting this work.}
\\Visilab Signal Technologies Oy, Finland}
\date{27th December, 2016}
\maketitle

\begin{abstract}
We transformed the generalized exponential power series to another functional form suitable for further analysis. By applying the Cauchy-Euler differential operator in the form of an exponential operator, the series became a sum of exponential differential operators acting on a simple exponential ($e^{-x}$). In the process we found new relations for the operator and a new polynomial with some interesting properties. Another form of the exponential power series became a nested sum of the new polynomial, thus isolating the main variable to a different functional dependence. We studied shortly the asymptotic behavior by using the dominant terms of the transformed series. New series expressions were created for common functions, like the trigonometric and exponential functions, in terms of the polynomial. \footnote{Visilab Report \#2016-12}
\subsection{Keywords}
exponential series, infinite series, Cauchy-Euler operator
\subsection{Mathematical Classification}
Mathematics Subject Classification 2010: 11L03, 30E10, 30K05, 30B10, 30C10, 30D10
\end{abstract}
\subsection{}

\tableofcontents{}
\section{Introduction}
\subsection{General}
\begin{figure}[ht]
	\centering
		\includegraphics[width=0.800\textwidth]{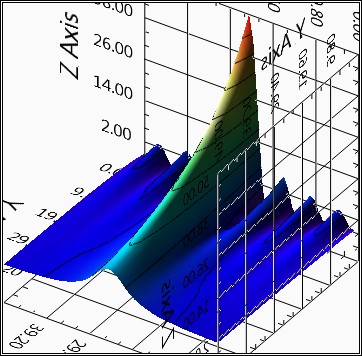}
	\caption{The real part of the generalized exponential series over the complex plane with $\alpha$=1}
	\label{fig:reseries1}
\end{figure}

\begin{figure}[ht]
	\centering
		\includegraphics[width=0.800\textwidth]{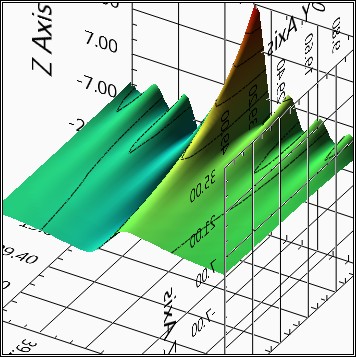}
	\caption{The imaginary part of the generalized exponential series over the complex plane with $\alpha$=1}
	\label{fig:imseries1}
\end{figure}

\begin{figure}[ht]
	\centering
		\includegraphics[width=0.800\textwidth]{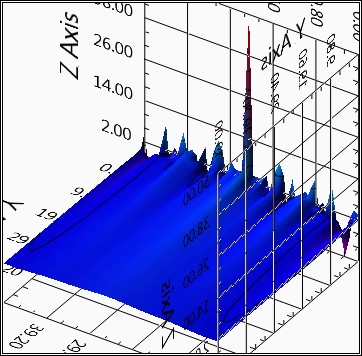}
	\caption{The real part of the generalized exponential series over the complex plane with $\alpha$=2.8}
	\label{fig:reseries28}
\end{figure}

\begin{figure}[ht]
	\centering
		\includegraphics[width=0.800\textwidth]{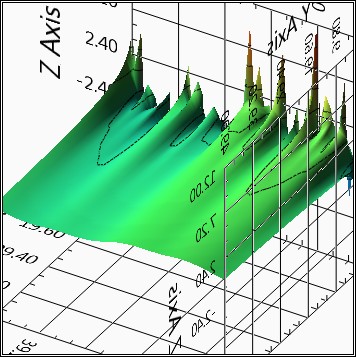}
	\caption{The imaginary part of the generalized exponential series over the complex plane with $\alpha$=2.8}
	\label{fig:imseries28}
\end{figure}

Our study focuses on a generalization of the quadratic exponential power series
\begin{equation}
f(x,2)=\sum^{\infty}_{k=1}e^{-x\cdot{k^2}} \label{eqn10}
\end{equation}
with $x>0,x\in{R}$. This function is very steep in behavior while $x$ is approaching zero, diverging fast and it goes quickly to zero as $x$ increases towards infinity. For negative $x$ values the series diverges. Our generalization will have a parameter $\alpha$
\begin{equation}
f(x,\alpha)=\sum^{\infty}_{k=1}e^{-x\cdot{k^{\alpha}}} \label{eqn20}
\end{equation} 
\begin{figure}[ht]
	\centering
		\includegraphics[width=1.00\textwidth]{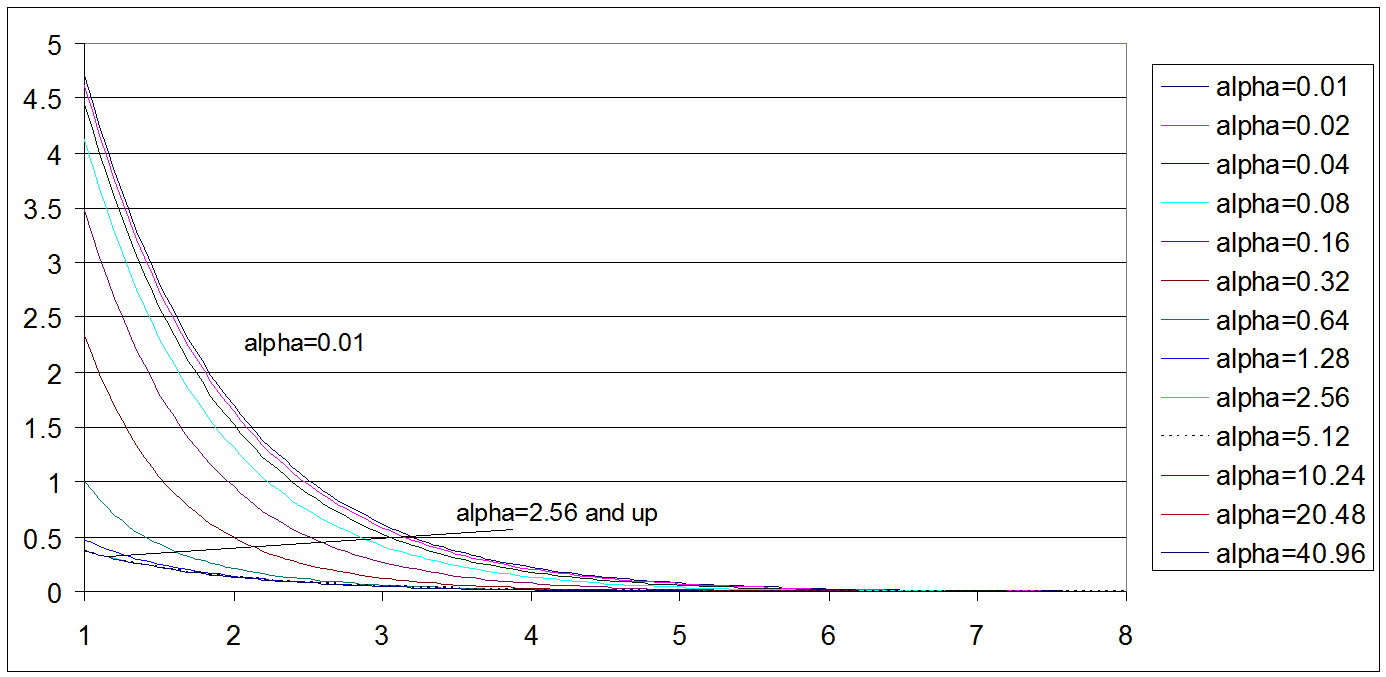}
	\caption{The generalized exponential series as a function of x with $\alpha$ as a parameter, excluding the origin, where it is singular}
	\label{fig:general}
\end{figure}
The range of validity is $x,\alpha\in{R}, x > 0$. By ratio test the generalized series converges absolutely when $\alpha > 0$. Figure \ref{fig:general} indicates how the family of curves saturates when $\alpha > 2.5$. The limit comes out as
\begin{equation}
\lim_{\alpha\rightarrow\infty} f(x,\alpha)=e^{-x}
\end{equation}
The generalized exponential series has a complicated behavior when the argument $x$ is varied over the complex plane. See the crude illustrations below (Figures \ref{fig:reseries1}, \ref{fig:imseries1}, \ref{fig:reseries28}, \ref{fig:imseries28}), generated with National Instruments LabVIEW2016. The origin is on the back plane lower corner center and the $x$ axis is directed towards the viewer. They show part of the behavior over first and fourth quadrants. We observe that the real part is behaving as an even function referring to the $x$ axis but the imaginary part is odd. The scales are only approximate. Rest of the illustrations are created with regular Excel spreadsheet graphics.

We are interested in studying the behavior within the ranges of $x,\alpha\in{R}$. We also study the asymptotic behavior at large $x\in{R}$. That is usually estimated as
\begin{equation}
f(x,\alpha){\approx}e^{-{x}}, x > 1 \label{eqn45}
\end{equation}
being the first term of the series. We wish to test if this approximation is superior to another estimate developed from the series. It is obvious that this is the main term at very large values of the argument. The questions are, what is the range of its validity and how good is our new approximation. In order to enter those tasks we need to transform or break down the exponential series in some manner.
\subsection{Attempts}
The exponential series effectively resists all common methods of summation to transform it to another form better suitable for analysis. Literature offers very little help in this respect. Jacobi theta functions are far relatives to this series. However, they are generally treated in different ways. We can offer a few na\'{i}ve attempts in order to penetrate into the internals. We can display the series term-wise and then collect them by taking a common factor at each step while progressing towards infinity.
\begin{equation}
f(x,2)=e^{-x}+e^{-4x}+e^{-9x}+e^{-16x}+e^{-25x}+e^{-36x}+...
\end{equation}
\begin{equation}
=e^{-x}(1+e^{-3x}(1+e^{-5x}(1+..)))
\end{equation}
Since with $n=0,1,2...$
\begin{equation}
(n+1)^2=\sum^{n}_{j=0}(2j+1) \label{eqn50}
\end{equation}
we get
\begin{equation}
f(x,2)=\sum^{\infty}_{k=1}{(e^{-x}e^{-3x}e^{-5x}e^{-7x}...e^{-(2k-1)x})} \label{eqn60}
\end{equation}
On the other hand we can first expand the series and then recollect to get
\begin{equation}
f(x,2)=\sum^{\infty}_{k=1}[{1-xk^2(1-\frac{xk^2}{2}(1-\frac{xk^2}{3}(1-\frac{xk^2}{4}(1-\frac{xk^2}{5}(1-\frac{xk^2}{6}(...))))))}] \label{eqn65}
\end{equation}
We can also force a regular Taylor's expansion to make it look like
\begin{equation}
f(x,2)=\sum^{\infty}_{k=1}{[1+\sum^{\infty}_{n=1}{\frac{(xk^2-2n)(xk^2)^{2n-1}}{(2n)!}}]} \label{eqn67}
\end{equation}
Unfortunately, these expressions do not seem to be helpful at this time.

The treatment of the subject starts in Chapter 2 by decomposing the original exponential series to nested series. Chapter 3 proceeds by transforming the series and a new polynomial is identified. It also offers recursion formulas and the generating function for it. Chapter 4 handles the differential operator acting on an exponential function which is equivalent to our exponential series. We handle there the eigenfunctions and eigenvalues of this operator. In Chapter 5 we study the asymptotic behavior of the series with the aid of the new results. Appendix \ref{apx:appa} displays the main properties of the Cauchy-Euler operator. In Appendix \ref{apx:appb} we study some properties of the new polynomial. Appendix \ref{apx:appc} treats shortly integration and differentiation properties of the series at hand. Appendix \ref{apx:appd} shows some results of generating new series expressions for some common functions. Our presentation is made in a  condensed way leaving out all formal proofs.
\section{Decomposition of the Series}
\subsection{Preliminaries}
We use the notation below while processing the expressions
\begin{equation}
z(k)=\alpha\cdot{ln(k)} \label{eqn1000}
\end{equation}
and
\begin{equation}
k^{\alpha}=e^{{\alpha}\cdot{ln(k)}}=e^z \label{eqn1010}
\end{equation}
Here $x,\alpha{\in{R}}$. The exponential nested structure needs to be broken down to some other functional dependence to allow easier handling in analysis.
\subsection{Exponential Series of an Exponential}
We start by expanding the exponential function as a Taylor's power series
\begin{equation}
f(x,\alpha)=\sum^{\infty}_{k=1}e^{-x\cdot{k^{\alpha}}}=\sum^{\infty}_{k=1, n=0}{\frac{(-x)^{n}\cdot{e^{nz}}}{n!}} \label{eqn1020}
\end{equation}
We expand further the term $e^{nz}$ as a power series getting
\begin{equation}
f(x,\alpha)=\sum^{\infty}_{k=1}{\sum^{\infty}_{n=0}{\sum^{\infty}_{m=0}{\frac{(-x)^{n}n^{m}z^{m}}{m!n!}}}} \label{eqn1030}
\end{equation}
Then we swap the summations $n,m$ to get
\begin{equation}
f(x,\alpha)=\sum^{\infty}_{k=1}{\sum^{\infty}_{m=0}{\sum^{\infty}_{n=0}{\frac{(-x)^{n}n^{m}z^{m}}{m!n!}}}} \label{eqn1040}
\end{equation}
This is also equal to 
\begin{equation}
f(x,\alpha)=\sum^{\infty}_{k=1}{[e^{-x}+\sum^{\infty}_{m=1}{\frac{z^{m}}{m!}\sum^{\infty}_{n=0}{\frac{(-x)^{n}n^{m}}{n!}}}]} \label{eqn1050}
\end{equation}
\subsection{Using a Function Instead of a Sum}
We define a function formed by the last nested sum as
\begin{equation}
S(x,m)=\sum^{\infty}_{n=0}{\frac{(-x)^{n}n^{m}}{n!}} \label{eqn1060}
\end{equation}
By using the property in Appendix \ref{apx:appa}. equation (\ref{eqn10020}), we can express it as
\begin{equation}
S(x,m)=(x\partial_x)^m\sum^{\infty}_{n=0}{\frac{(-x)^{n}}{n!}}=(x\partial_x)^{m}e^{-x} \label{eqn1070}
\end{equation}
Then we can write our function as
\begin{equation}
f(x,\alpha)=\sum^{\infty}_{k=1}{\sum^{\infty}_{m=0}{\frac{z^{m}}{m!}S(x,m)}} \label{eqn1080}
\end{equation}
and equivalently as follows
\begin{equation}
f(x,\alpha)=\sum^{\infty}_{k=1}{\sum^{\infty}_{m=0}{\frac{z^{m}}{m!}(x\partial_x)^{m}e^{-x}}} \label{eqn1090}
\end{equation}
\begin{equation}
=\sum^{\infty}_{k=1}{K(x,z)}e^{-x} \label{eqn1100}
\end{equation}
$K(x,z)$ is a differential operator (recall that $z(k)=\alpha\cdot{ln(k)}$)
\begin{equation}
K(x,z)=\sum^{\infty}_{m=0}{\frac{(zx\partial_x)^{m}}{m!}}=e^{zx\partial_x} \label{eqn1110}
\end{equation}
\section{Polynomials}
\subsection{Triangle}
\begin{figure}[ht]
	\centering
		\includegraphics[width=1.00\textwidth]{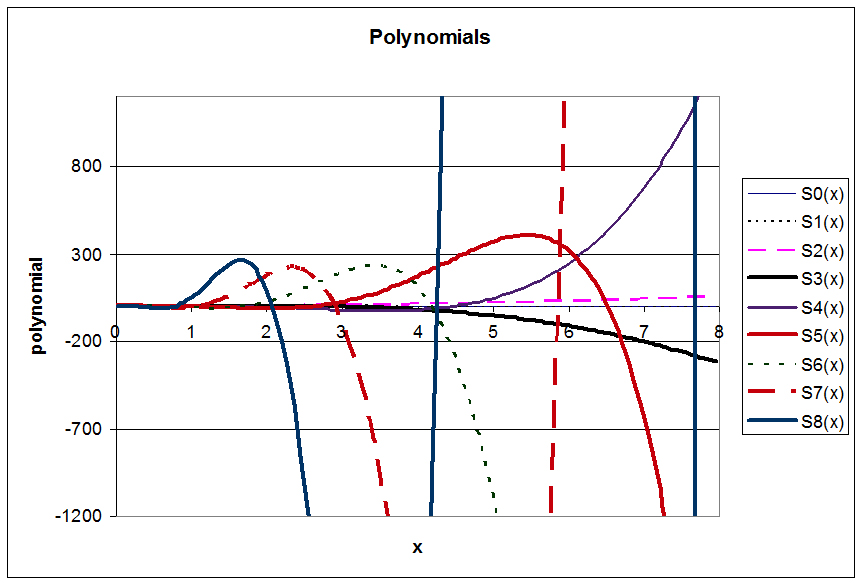}
	\caption{Polynomials $S_m(x)$}
	\label{fig:poly}
\end{figure}

\begin{figure}[ht]
	\centering
		\includegraphics[width=1.00\textwidth]{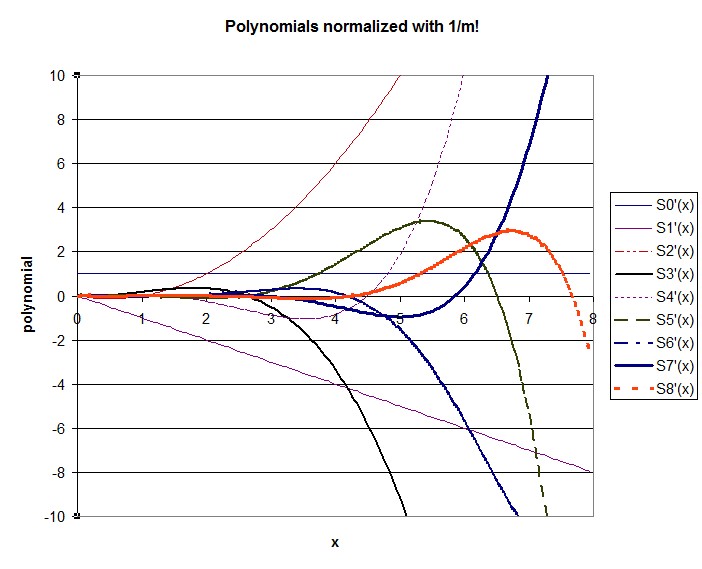}
	\caption{Normalized polynomials as $\frac{S_m(x)}{m!}$}
	\label{fig:polynorm}
\end{figure}

\begin{figure}[ht]
	\centering
		\includegraphics[width=1.00\textwidth]{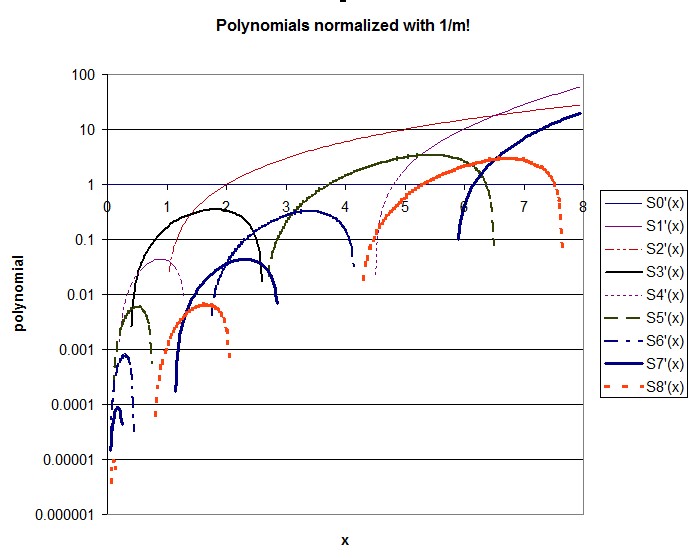}
	\caption{Normalized polynomials as $\frac{S_m(x)}{m!}$ in log scale. Note that the logarithmic illustration has cut the negative values.}
	\label{fig:polynormlog}
\end{figure}

We can multiply the function $S(x,m)$ in equation (\ref{eqn1070}) by $e^x$ from the left and mark it as $S_m(x)$
\begin{equation}
S_m(x)=e^{x}S(x,m)=e^{x}(x\partial_x)^{m}e^{-x} \label{eqn2000}
\end{equation}
Thus the exponential functions disappear for every $m$ and $S_m(x)$ appears to be a new kind of polynomial. In the following we study a little of its properties. It is notable that $S_m(x)$ and $S(x,m)$ are independent of $\alpha$. From equation (\ref{eqn2000}) we can solve for the first polynomials by direct differentiation. We get the triangle below.
\begin{equation}
S_0(x)=1 \nonumber
\end{equation}
\begin{equation}
S_1(x)=-x \nonumber
\end{equation}
\begin{equation}
S_2(x)=x^2-x \nonumber
\end{equation}
\begin{equation}
S_3(x)=-x^3+3x^2-x \nonumber
\end{equation}
\begin{equation}
S_4(x)=x^4-6x^3+7x^2-x \nonumber
\end{equation}
\begin{equation}
S_5(x)=-x^5+10x^4-25x^3+15x^2-x \nonumber
\end{equation}
\begin{equation}
S_6(x)=x^6-15x^5+65x^4-90x^3+31x^2-x \nonumber
\end{equation}
\begin{equation}
S_7(x)=-x^7+21x^6-140x^5+350x^4-301x^3+63x^2-x \nonumber
\end{equation}
\begin{equation}
S_8(x)=x^8-28x^7+266x^6-1050x^5+1701x^4-966x^3+127x^2-x \label{eqn2010}
\end{equation}
\begin{equation}
... \nonumber
\end{equation}
Figures \ref{fig:poly}, \ref{fig:polynorm} and \ref{fig:polynormlog} illustrate the first few polynomials, the latter normalized with $m!$. The polynomial coefficients resemble the Stirling's numbers of second kind but with a changing sign.

\subsection{Recursion Relations for the Polynomials}
We can find a recursion relation for the coefficients of the powers $b_{j}^{n}$, ($j$ row, $n$ power), by looking at the elements on the row immediately above in the triangle
\begin{equation}
(-1)^{n+j}b_{j}^{n}x^{n-j}+(-1)^{n+j+1}b_{j+1}^{n}x^{n-j-1}+...  \label{eqn2030}
\end{equation}
\begin{equation}
...(-1)^{n+j+1}b_{j+1}^{n+1}x^{n+1-j-1}+...  \label{eqn2040}
\end{equation}
We get
\begin{equation}
b_{j+1}^{n+1}=(n-j){b_{j}^{n}}+{b_{j+1}^{n}}  \label{eqn2050}
\end{equation}
where we have the absolute value of the coefficient only. This recursion is too clumsy for use in analysis. In the following we develop proper recursion relations for the $S_m(x)$. The simplest one follows directly from equation (\ref{eqn2000}) by multiplying it from the left with $x\partial_x$.
\begin{equation}
xS_n^{'}(x)=xS_n(x)+S_{n+1}(x) \ \ \ \ \ \ (RR1) \label{eqn2060}
\end{equation}
or
\begin{equation}
S_{n+1}(x)=x(S_n^{'}(x)-S_n(x)) \label{eqn2070}
\end{equation}
\subsection{Generating Function of the Polynomial}
An important step is to find the generating function for the polynomial $S_m(x)$. By multiplying equation (\ref{eqn2000}) from the left by $e^{-x}$, we get
\begin{equation}
S_{m}(x)e^{-x}=(x\partial_x)^{m}e^{-x} \label{eqn2080}
\end{equation}
Next we multiply both sides with $\frac{t^m}{m!}$ and sum the terms to obtain
\begin{equation}
\sum^{\infty}_{m=0}{\frac{S_{m}(x)e^{-x}t^{m}}{m!}}=\sum^{\infty}_{m=0}{\frac{(x\partial_x)^{m}e^{-x}t^{m}}{m!}}=e^{tx\partial{x}}e^{-x} \label{eqn2090}
\end{equation}
leading to the generating function
\begin{equation}
\sum^{\infty}_{m=0}{\frac{S_{m}(x)e^{-x}t^{m}}{m!}}=e^{-xe^{t}} \ \ \ \ \ \ (GF) \label{eqn2100}
\end{equation}
We have used the Cauchy-Euler operator properties as expressed in Appendix \ref{apx:appa}.
\subsection{Making Symmetric the Generating Function}
It turns out that we can make symmetric the generating function (GF) by changing the variable $t$ to $y$ as
\begin{equation}
t=ln(y+1)
\end{equation}
we are able to write down the generating function in the following form
\begin{equation}
e^{-xy}=\sum^{\infty}_{m=0}{\frac{S_{m}(x)(ln(y+1))^{m}}{m!}} \label{eqn2102}
\end{equation}
The symmetry is equal to the fact that on the left side we can swap $y\rightleftharpoons{x}$. Therefore, we can do the same swapping on the right side without changing the value of the expression. This is a rare symmetry property. In general, a function 
\begin{equation}
f(xy)=g(x,y)
\end{equation}
can be partially differentiated to obtain
\begin{equation}
x\partial{x}g(x,y)=y\partial{y}g(x,y)
\end{equation}
thus containing the Cauchy-Euler differential operator.
\subsection{Further Recursion Relations}
By defining
\begin{equation}
S_m(x)=\sum^{m}_{j=1}{c_{j}^{m}x^{m-j+1}}  \label{eqn2200}
\end{equation}
we get by using the equation (\ref{eqn2070})
\begin{equation}
c_j^m=c_{j-1}^{m-1}(m-j+1)-c_{j}^{m-1}  \label{eqn2210}
\end{equation}
The values below work as boundaries
\begin{equation}
c_m^m=-1  \label{eqn2230}
\end{equation}
\begin{equation}
c_1^m=(-1)^m  \label{eqn2240}
\end{equation}
\begin{equation}
c_j^m=0, j>m, j<1  \label{eqn2250}
\end{equation}
Even this recursion relation in unsatisfactory. We will obtain a pair of good recursion relations by differentiating the (GF) equation (\ref{eqn2100}) with $\partial_x$ and $\partial_t$ separately.
\begin{equation}
-x{e^t}\sum^{\infty}_{m=0}{\frac{S_{m}(x)e^{-x}t^{m}}{m!}}=\sum^{\infty}_{m=0}{\frac{S_{m+1}(x)e^{-x}t^{m}}{m!}}  \label{eqn2252}
\end{equation}
\begin{equation}
-{e^t}\sum^{\infty}_{m=0}{\frac{S_{m}(x)e^{-x}t^{m}}{m!}}=\sum^{\infty}_{m=0}{\frac{(S_{m}^{'}(x)-S_{m}(x))e^{-x}t^{m}}{m!}}  \label{eqn2254}
\end{equation}
Comparison of the powers of $t$ brings out
\begin{equation}
\sum^{j}_{n=0}{\frac{S_{j-n}(x)}{(j-n)!n!}}=\frac{S_{j}(x)-S_{j}^{'}(x)}{j!}  \ \ \ \ \ \ (RR2) \label{eqn2260}
\end{equation}
and 
\begin{equation}
-x\sum^{j}_{n=0}{\frac{S_{j-n}(x)}{(j-n)!n!}}=\frac{S_{j+1}(x)}{j!}  \ \ \ \ \ \ (RR3) \label{eqn2270}
\end{equation}
The latter recursion is useful for generation of the polynomials.
\subsection{Explicit Form of the Polynomial}
By looking at the recursion relations obtained, we recognize that the Stirling number of the second kind $\hat{S}_n^{[m]}$ is reminiscent to our polynomial coefficients. It has a closed-form expression as follows, (\cite{Gradshteyn2007}, \cite{Abramowitz1970})
\begin{equation}
\hat{S}_n^{[m]}=\frac{1}{m!}\sum^{m}_{i=0}{(-1)^{m-i}(^m_i)i^n} \ \ \ \ \  \label{eqn2275}
\end{equation}
and a recursion formula
\begin{equation}
\hat{S}_{n+1}^{[m]}=m\hat{S}_n^{[m]}+\hat{S}_n^{[m-1]}  \label{eqn2277}
\end{equation}
It follows that we are able to express our polynomial in terms of it
\begin{equation}
S_m(x)=\sum^{m-1}_{j=0}{(-x)^{m-j}\hat{S}^{[m-j]}_m} \ \ \ \ \  \label{eqn2279}
\end{equation}
Thus we get an explicit form
\begin{equation}
S_m(x)=\sum^{m-1}_{j=0}{(-x)^{m-j}\frac{1}{(m-j)!}\sum^{m-j}_{i=0}{(-1)^{m-j-i}(^{m-j}_i)i^m}} \ \ \ \ \  \label{eqn2300}
\end{equation}
shortening to the final expression for the polynomial
\begin{equation}
S_m(x)=\sum^{m-1}_{j=0}{x^{m-j}\sum^{m-j}_{i=0}{\frac{(-1)^{i}i^m}{i!(m-j-i)!}}},   m=1,2,3.. \ \ \ \  \label{eqn2310}
\end{equation}
\begin{equation}
S_0(x)=1
\end{equation}
\section{Operator Expressions}
\subsection{The Differential Operator}
To collect our results thus far, we obtain (recalling that $z=z(k)$) and by using equation (\ref{eqn1110})
\begin{equation}
f(x,\alpha)=\sum^{\infty}_{k=1}{e^{zx\partial_x}e^{-x}} \label{eqn2710}
\end{equation}
\begin{equation}
=\sum^{\infty}_{k=1}{e^{\alpha{ln(k)}x\partial_x}e^{-x}} \label{eqn2720}
\end{equation}
\begin{equation}
=\sum^{\infty}_{k=1}{\frac{1}{k^{-\alpha{x}\partial_x}}e^{-x}} \label{eqn2730}
\end{equation}
\begin{equation}
=\zeta(-\alpha{x}\partial_x)e^{-x} \label{eqn2740}
\end{equation}
$\zeta(x)$ is formally the Riemann zeta function, now becoming a differential operator. Since the expansion can be made inside the k-summation (\ref{eqn1020}), we have finally everything collected to our transformed series
\begin{equation}
f(x,\alpha)=\sum^{\infty}_{k=1}{e^{-xk^{\alpha}}}=\sum^{\infty}_{k=1}{k^{\alpha{x}\partial_x}e^{-x}}=e^{-x}\sum^{\infty}_{k=1}{\sum^{\infty}_{m=0}{\frac{(\alpha{ln(k)})^m}{m!}S_{m}(x)}} \label{eqn2760}
\end{equation}
Thus we have managed to move the $x$-dependence to the polynomial and also the $\alpha$-dependence is isolated from that. On the other hand, we can apply the (GF) to the inner sum to write down
\begin{equation}
e^{-x}\sum^{\infty}_{k=1}{\sum^{\infty}_{m=0}{\frac{(\alpha{ln(k)})^m}{m!}S_{m}(x)}}=e^{-x}\sum^{\infty}_{k=1}{e^{-xe^{\alpha{ln(k)}}+x}} 
\end{equation}
This becomes immediately an identity as expected. This once more proves that the generating function is the essential core for our generalized exponential series.

\subsection{Eigenvalues and Eigenfunctions of the Differential Operator}
The eigenfunctions and eigenvalues of the operator in (\ref{eqn2740}) are determined as follows.
\begin{equation}
\zeta(-\alpha{x}\partial_x)\phi{(x)}=\sum^{\infty}_{k=1}{\phi{(xe^{\alpha{ln(k)}})}}=\sum^{\infty}_{k=1}{\phi{(xk^{\alpha}}}) \label{eqn2800}
\end{equation}
As a trial function we set, with $\beta{>0}$, a constant. $Re{(\alpha{\beta})}>1$
\begin{equation}
\zeta(-\alpha{x}\partial_x)\frac{1}{x^{\beta}}=\frac{1}{x^{\beta}}\sum^{\infty}_{k=1}{k^{-\alpha{\beta}}}=\frac{\zeta(\alpha{\beta})}{x^{\beta}} \label{eqn2820}
\end{equation}
This is convergent if $Re{(\alpha{\beta})}>1$ and we have the essential result of a continuum of eigenvalues and eigenfunctions for the $\zeta()$ operator. The $N$ eigenvalues are thus $\zeta(\alpha{\beta})$ and eigenfunctions are $\frac{1}{x^{\beta}}$. 

If $\phi(x)$ can be expanded as a kind of a Laurent series with negative powers only, we have, with $\phi{_j}$ constant coefficients of the expansion
\begin{equation}
\phi(x)=\sum^{N}_{j=1}{\frac{\phi_{j}}{x^j}} \label{eqn2830}
\end{equation}
and therefore, by (\ref{eqn2820})
\begin{equation}
\zeta(-\alpha{x}\partial_x)\phi(x)=\sum^{N}_{j=1}{\frac{\phi_{j}\zeta(\alpha{j})}{x^j}} \label{eqn2840}
\end{equation}
\section{Asymptotic Behavior}
\subsection{Large Argument Estimate}
\begin{figure}[ht]
	\centering
		\includegraphics[width=1.00\textwidth]{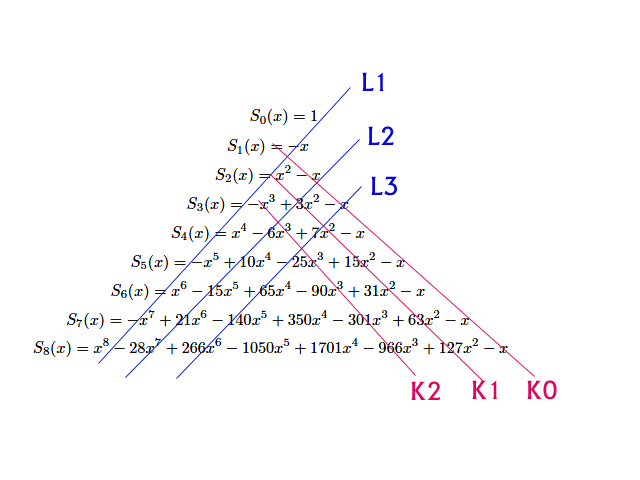}
	\caption{Definitions of the diagonals, left and right slanted, for the polynomial triangle}
	\label{fig:diagonals}
\end{figure}

We can apply some of the results obtained to estimate the asymptotic behavior of the exponential series equation (\ref{eqn20}) when the argument $x$ approaches large values. We can look at the polynomial triangle equation (\ref{eqn2010}) along diagonals, both slanted to the left and to the right, Figure \ref{fig:diagonals}. We mark the right-slanted diagonals as $L_n$ and left-slanted diagonals as $K_n$. The first two can be identified as closed-form expressions but the diagonals following them will have very complicated expressions. Equation (\ref{eqn2310}) can be used directly as well if more terms are required. The terms are, while $m$ is referring to the power of the polynomial, the in first right-slanted diagonal 
\begin{equation}
L_1(m)= (-1)^{m}x^m, \ \ \ \ \ m=0,1,2,3... \label{eqn4000}
\end{equation}
and in the second one
\begin{equation}
L_2(m)= (-1)^{m+1}x^{m-1}m(m-1)\frac{1}{2},  \ \ \ \ \  m=2,3... \label{eqn4010}
\end{equation}
The higher coefficients are likely too complicated to be evaluated just by looking at the triangle. The general coefficients can be obtained from equation (\ref{eqn2310}) with $j=n-1$
\begin{equation}
L_n(m)=\sum^{m-n+1}_{i=0}{\frac{(-1)^{i}i^{m}x^{m-n+1}}{i!(m-n+1-i)!}}
\end{equation}

The largest terms are those with the highest powers and we pick up $L_1$ and $L_2$ terms for the $S_m(x)$. The exponential series will become
\begin{equation}
f(x,\alpha)=\sum^{\infty}_{k=1}{e^{-xk^{\alpha}}}=\sum^{\infty}_{k=1}{\sum^{\infty}_{m=0}{\frac{(\alpha{ln(k)})^m}{m!}S_{m}(x)e^{-x}}} \label{eqn4100}
\end{equation}

\begin{equation}
\approx{e^{-x}\sum^{\infty}_{k=1}{[\sum^{\infty}_{m=0}{\frac{(\alpha{ln(k)})^{m}(-x)^m}{m!}}-\frac{1}{x}\sum^{\infty}_{m=2}{\frac{(-\alpha{{x}ln(k)})^{m}m(m-1)}{2(m!)}}]}} \label{eqn4110}
\end{equation}
Instantly we recognize the familiar Riemann zeta function and get
\begin{equation}
f(x,\alpha)\approx{e^{-x}\zeta(x\alpha)-\frac{e^{-x}x}{2}\zeta^{''}{(x\alpha)}}, x >> 0 \label{eqn4120}
\end{equation}
The range of validity can be estimated below with a graph. The traditional estimate has been the first term only
\begin{equation}
f(x,\alpha){\approx}e^{-{x}} \label{eqn4200}
\end{equation}
We compare this to the actual function and to the estimate (\ref{eqn4120}) in the following Figures \ref{fig:asymptotics} and \ref{fig:asymptoticslog}.
\begin{figure}[ht]
	\centering
		\includegraphics[width=1.0\textwidth]{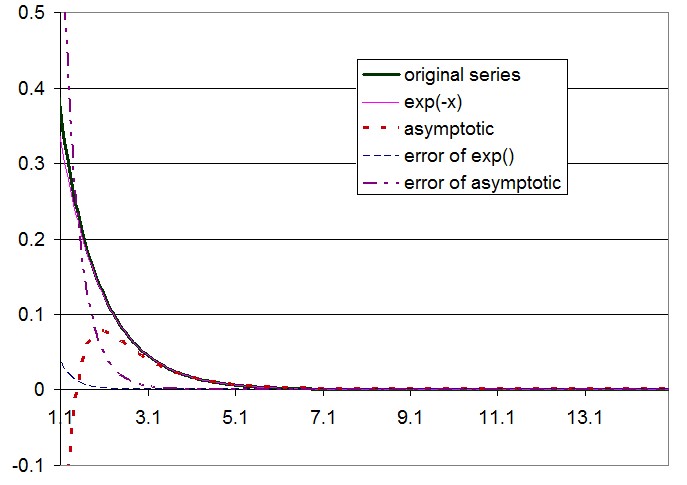}
	\caption{Asymptotic behavior of the exp(-x) estimate and our approximation of two terms, with $\alpha$=1.6, for large argument values, linear scale}
	\label{fig:asymptotics}
\end{figure}
\begin{figure}[ht]
	\centering
		\includegraphics[width=1.0\textwidth]{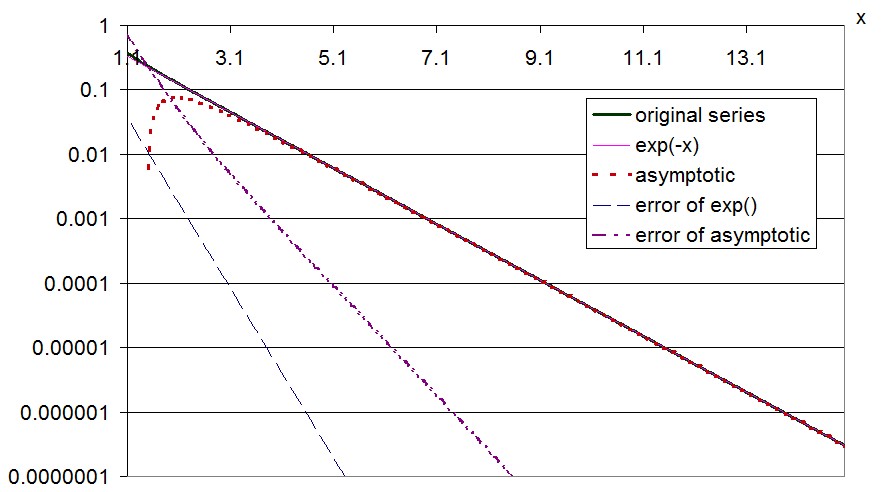}
	\caption{Asymptotic behavior of the exp(-x) estimate and our approximation of two terms, with $\alpha$=1.6, for large argument values, logarithmic vertical scale}
	\label{fig:asymptoticslog}
\end{figure}

The graphs show that for large values of $x > 1$ the simple exponential estimate still gives a better accuracy than our asymptotic model with two biggest terms. We would need more terms for the approximation to possibly beat the simple exponential function. However, that is not justified since the new terms will increase complexity too much. 

\newpage
\subsection{Asymptotic Behavior Near Zero}
To select the lowest powers in $x$ we pick up the left-slanted diagonals $K_n(m)$ (see Figure \ref{fig:diagonals})
\begin{equation}
K_0(m)= -x,   \ \ \ \ \  m=1,2... \label{eqn4710}
\end{equation}
\begin{equation}
K_1(m)= (2^{m-1}-1)x^2,  \ \ \ \ \  m=2,3... \label{eqn4720}
\end{equation}
Again, higher terms become awkward to evaluate by guessing from the triangle and the general expression is
\begin{equation}
K_n(m)=\sum^{n+1}_{i=0}{\frac{(-1)^{i}i^{m}x^{m+1}}{i!(n+1-i)!}}
\end{equation}
Thus
\begin{equation}
K_2(m)= (2^{m}-1-3^{m-1})\frac{x^3}{2!},  \ \ \ \ \  m=3,4... \label{eqn4722}
\end{equation}
\begin{equation}
K_3(m)= (3\cdot{2^{m-1}}-1-3^m+4^{m-1})\frac{x^4}{3!},  \ \ \ \ \  m=4,5... \label{eqn4724}
\end{equation}
$K_n(m)$ refers always to the same power $m$ of $x$.
In order to estimate the behavior when $x$ approaches zero, we may take the terms $K_0$ and $K_1$ to our expressions since they are with the lowest powers. This case is much more complicated than the preceding one as the series will always start to diverge while approaching the limit. The question is only, how does it do it?
The exponential series is
\begin{equation}
f(x,\alpha)=\sum^{\infty}_{k=1}{e^{-xk^{\alpha}}}=e^{-x}\sum^{\infty}_{k=1}{\sum^{\infty}_{m=0}{\frac{(\alpha{ln(k)})^m}{m!}S_{m}(x)}} \label{eqn4770}
\end{equation}

\begin{equation}
\approx{e^{-x}(\sum^{\infty}_{k=1}{1+\frac{x^2}{2}+x-(x^2+x)k^{\alpha}+\frac{x^2k^{2\alpha}}{2})}} \label{eqn4780}
\end{equation}
The sum will diverge rapidly while $x$ approaches zero. The main term is the unity and the other terms will be big too but of opposite signs. At this end of the range, we have much more trouble in obtaining a simple approximating function for the generalized exponential function. Since this approach fails, the analysis must be carried out in another way to obtain an asymptotic function near zero. 

We are able to investigate the behavior at small $x$ in the case of $\alpha=1$ since we can solve the simple series as follows
\begin{equation}
\sum^{\infty}_{k=1}{e^{-kx}}=\frac{e^{-x}}{1-e^{-x}} 
\end{equation}
Taking the limit of $x$ approaching zero shows that the result is a simple pole. We are not with the liberty of claiming that this would be valid for other values of $\alpha$.

\section{Conclusions}
The exponential series and its generalized version having a power of the index $\alpha$, appeared to be a hard nut to be opened. Our aim was to transform the generalized exponential series, to a form with a simpler functional dependence of the parameter $x$. 

We first decomposed the series by making it a series of nested exponents and restructured it into nested series instead. We recognized the new form by applying an exponential Cauchy-Euler differential operator. This allowed further transformation of the innermost power series to become a polynomial. This could, on the other hand, be understood as a series of powers of a differential operator acting on an exponential function. This series was formally identified as the Riemann zeta function with an argument $-\alpha{x}\partial{x}$, simplifying presentation. The original series was converted into a nested series of polynomials in $x$, allowing a simplification for further analysis.

While developing the transformation, we found a new polynomial with interesting properties. Those are expressed in the Appendices and some results are also in the main text. The recursion relations for the polynomial were solved, equations (\ref{eqn2060}), (\ref{eqn2260}) and (\ref{eqn2270}). With the aid of them, we presented a self-contained expression for the polynomial, equation (\ref{eqn2310}) which can be used for further analysis of the generalized exponential series. The polynomial can be simply solved by hand with the aid of a triangle as in equations (\ref{eqn2000}), (\ref{eqn2010}) and (\ref{eqn2050}). 

The generating function for the polynomial is a central tool for unwrapping the properties of the polynomial, equation (\ref{eqn2100}). It is based on the first recursion relation (\ref{eqn2060}). Some new series expressions for the common exponential function and for the two basic trigonometric functions are in Appendix \ref{apx:appd}. They were derived with the aid of the generating function. We made a change of variable to the generating function forcing it to become symmetric with respect to the two parameters, equation (\ref{eqn2102}). This is a series of two multiplying functions where the two parameters are separated and can be swapped on both sides. This led to some new interesting series expressions shown in Appendix \ref{apx:appd}.

One of our motives was to better understand the asymptotic behavior at large variable values and also near zero. We tested whether the transformed series can offer a satisfactory approximation. It appeared that for large $x$ by using a few dominant terms of the new series, it fails to give a better asymptotic estimate than the traditional single-term approximation ($=e^{-x}$). While approaching zero from the positive side, there seems not to be any available approximation from our analysis. More work is required in this field.

\appendix
\section{Appendix. Properties of the Exponential Cauchy-Euler Operator}\label{apx:appa}
Assuming $x\in{R}$, we present some useful features of the Cauchy-Euler differential operator. It is applied to functions whose derivatives exist and we also assume that the functions can be expanded as a Taylor's power series (or any other power series with positive integer powers). The operator appears to have interesting properties displayed here, most of which are believed to be new. Here $\partial_{x}$ represents the partial derivative operator. The results are valid for analytic functions over the complex plane as well.
\subsection{Basic Operations}
We begin by presenting progressively increasing complexity for the operator. These are proven by differentiation and are well known.
\begin{equation}
(x\partial_{x})x=x \label{eqn9000}
\end{equation}
\begin{equation}
(x\partial_{x})^m{x}=x \label{eqn9010}
\end{equation}
\begin{equation}
(x\partial_{x})x^n=n\cdot{x^n} \label{eqn10000}
\end{equation}
Therefore we obtain
\begin{equation}
(x\partial_{x})^m{x^n}=n^m\cdot{x^n} \label{eqn10020}
\end{equation}
For the general case, $m\in{N}$, $n^{m}$ is the eigenvalue and $x^n$ is the eigenfunction, $n\in{R}, n\neq{0}$.

The partial derivative and Cauchy-Euler operator have commutators which are useful in handling the Cauchy-Euler operator in more complex cases.
\begin{equation}
[x,\partial_{x}]=-1 \label{eqn10025}
\end{equation}
\begin{equation}
[x\partial_{x},\partial_{x}{x}]=0 \label{eqn10030}
\end{equation}
We apply the $m$'th power operator to a more challenging function and get after expanding the exponential function to a Taylor's power series
\begin{equation}
(x\partial_{x})^m{e^{-x}}=\sum^{\infty}_{n=0}{\frac{(-x)^{n}n^{m}}{n!}} \label{eqn10040}
\end{equation}
We wish to expand the exponential Cauchy-Euler differential operator and define it by its power series. That leads to
\begin{equation}
e^{\beta{x}\partial_{x}}{x^j}=\sum^{\infty}_{n=0}{\frac{\beta^{n}(x\partial_{x})^n}{n!}}{x^{j}}=\sum^{\infty}_{n=0}{\frac{\beta^{n}j^n{x^{j}}}{n!}} \label{eqn10050}
\end{equation}
\begin{equation}
=e^{\beta{j}}x^{j}=(e^{\beta}x)^{j} \label{eqn10060}
\end{equation}
Here we have taken into use a parameter $\beta{\in}C$. 
\subsection{Subjected to a General Function}
Without losing generality, we assume the function $A(x)$ has a Taylor's series around the origin. Thus, by using the results above, we can write down a very useful general expression
\begin{equation}
e^{\beta{x}\partial_{x}}A(x)=e^{\beta{x}\partial_{x}}\sum^{\infty}_{n=0}{\frac{a_n{x}^n}{n!}}=\sum^{\infty}_{n=0}{\frac{{a_n}e^{\beta{x}\partial_{x}}{x}^{n}}{n!}} \nonumber
\end{equation}
\begin{equation}
=\sum^{\infty}_{n=0}{\frac{{a_n}{(e^{\beta}x)^{n}}}{n!}}=A(xe^\beta) \label{eqn10070}
\end{equation}
We have as simple examples of application
\begin{equation}
cos({\beta{x}\partial_{x}})A(x)=\frac{1}{2}[A(xe^{i\beta})+A(xe^{-i\beta})] \label{eqn10100}
\end{equation}
\begin{equation}
sin({\beta{x}\partial_{x}})A(x)=\frac{1}{2i}[A(xe^{i\beta})-A(xe^{-i\beta})] \label{eqn10110}
\end{equation}
In particular we have
\begin{equation}
e^{\beta{x\partial_{x}}}e^{-\eta{x}}=e^{-\eta{xe^{\beta}}} \label{eqn10180}
\end{equation}
The equation (\ref{eqn10070}) will be useful while transforming the generalized exponential series. It may also find use in other fields, like in solving the Cauchy problem and in transforming or linearizing differential equations.

The inverse of the exponential Cauchy-Euler operator can be found as
\begin{equation}
e^{-\alpha{x}\partial_x}(e^{\alpha{x}\partial_x}A(x))=A(x) 
\end{equation}
or
\begin{equation}
e^{-\alpha{x}\partial_x}(e^{\alpha{x}\partial_x})=I 
\end{equation}
which is expected. 

For comparison, we write below the known results of the ordinary exponential differential operator. $a $ is a constant.
\begin{equation}
e^{a\partial_{x}}x^n=(x+a)^n
\end{equation}
\begin{equation}
e^{a\partial_{x}}A(x)=A(x+a) 
\end{equation}

\subsection{Subjected to a Polynomial by a General Operator Function}
If a function $\phi(x)$ can be expanded as a finite or infinite power series around origin
\begin{equation}
\phi(x)=\sum^{M}_{n=0}{x^n{f_n}} \label{eqn14880}
\end{equation}
with $f_n$ constant coefficients, we may be able to have an expression with an operator function having as well a basic finite or infinite polynomial expansion, with $b_j$ constants and no $x$ dependence.
\begin{equation}
b(u)=\sum^{L}_{j=0}{u^j{b_j}} \label{eqn14890}
\end{equation}
as follows
\begin{equation}
b(\beta{x}\partial_x)\phi(x)=\sum^{M}_{n=0}{x^n{f_n}b(\beta{n})} \label{eqn14900}
\end{equation}
\subsection{Subjected to an Arbitrary Function by a General Operator Function}
If we can define an operator function $B(x)$ as a regular Taylor's series with $b_n$ as constant coefficients, we may express it as
\begin{equation}
B(u)=\sum^{\infty}_{n=0}{u^n{b_n}} \label{eqn15890}
\end{equation}
This operator can now be applied to the function $A(x)$ (having a series expansion as above) with a result as follows
\begin{equation}
B(e^{\beta{x}\partial{x}})A(x)=\sum^{\infty}_{n=0}{{b_n}A(xe^{n\beta})} \label{eqn15880}
\end{equation}
\section{Appendix. Properties of the New Polynomial}\label{apx:appb}
\subsection{Function with the Exponential Factor}
We have as the starting point for the function containing the exponential factor,
\begin{equation}
S(x,m)=(x\partial_x)^{m}e^{-x} \label{eqn11000}
\end{equation}
we can prove by integration
\begin{equation}
\int_0^{\infty}{dx{S(x,0)}}=1  \label{eqn11010}
\end{equation}
\begin{equation}
\int_0^{\infty}{dx{S(x,1)}}=-1  \label{eqn11020}
\end{equation}
\begin{equation}
\int_0^{\infty}{dx{S(x,2)}}=1  \label{eqn11030}
\end{equation}
or in general
\begin{equation}
\int_0^{\infty}{dx{S(x,m)}}=(-1)^m  \label{eqn11040}
\end{equation}
\subsection{Polynomial Identities}
By using the earlier summation formula (\ref{eqn3890}), we get after differentiating this by $\partial_x$
\begin{equation}
\sum^{\infty}_{m=0}{\frac{S^{'}_{m}(x)}{m!}}=-(e-1)e^{-x(e-1)}  \label{eqn11060}
\end{equation}
We can apply the recursion relation (RR1) to this and get
\begin{equation}
\sum^{\infty}_{m=0}{\frac{(S_{m+1}(x)+exS_m(x))}{m!}}=0  \label{eqn11070}
\end{equation}
This is equal to 
\begin{equation}
\sum^{\infty}_{m=0}{\frac{S_{m+1}(x)}{m!}}=-ex{e^{-x(e-1)}}  \label{eqn11080}
\end{equation}
By differentiating the equation (\ref{eqn11070}) and by using the (RR1) again we get
\begin{equation}
\sum^{\infty}_{m=0}{\frac{((1+e)S_{m+1}(x)+(1+x)eS_m(x)+\frac{1}{x}S_{m+2}(x))}{m!}}=0  \label{eqn11090}
\end{equation}
We can place equation (\ref{eqn11070}) into this to get
\begin{equation}
\sum^{\infty}_{m=0}{\frac{(eS_{m+1}(x)+eS_m(x)+\frac{1}{x}S_{m+2}(x))}{m!}}=0  \label{eqn11100}
\end{equation}
Placing it again will produce the following
\begin{equation}
\sum^{\infty}_{m=0}{\frac{((1-ex)S_{m}(x)+\frac{1}{ex}S_{m+2}(x))}{m!}}=0  \label{eqn11110}
\end{equation}
We can continue this kind of processing to get
\begin{equation}
\sum^{\infty}_{m=0}{\frac{S_{m+3}(x)}{m!}}=ex(3ex-1-e^2{x^2})\sum^{\infty}_{m=0}{\frac{S_{m}(x)}{m!}}  \label{eqn11120}
\end{equation}
\begin{equation}
=ex(3ex-1-e^2{x^2})e^{-x(e-1)}  \nonumber
\end{equation}
By differentiating with $\partial_t$ the generating  function (\ref{eqn2100}), we obtain after setting $t=1$
\begin{equation}
\sum^{\infty}_{m=0}{\frac{S_{m+2}(x)}{m!}}=(-xe+x^2e^2)e^{-x(e-1)}  \label{eqn11160}
\end{equation}
Obviously, there is no end to generating new equations in this way. 

The last sum below is rather trivial.
\begin{equation}
\sum_{j=0}^{m-1}(-1)^{j+1}(S_j+S_{j+1})=-1+(-1)^{m}S_m 
\end{equation}
\subsection{Polynomial Integrated}
Partial integration will give a useful result
\begin{equation}
\int{dx{e^{-x}S_{m}(x)}}=(-1)^{m}[xe^{-x}\sum^{m-1}_{j=0}{(-1)^{j+1}S_j(x)}-e^{-x}]\label{eqn13410}
\end{equation}
and (RR1) gives
\begin{equation}
\int{dx{e^{-x}S_{m+1}(x)}}=xe^{-x}S_m(x)-\int{{dx}e^{-x}S_m(x)}\label{eqn13420}
\end{equation}
\subsection{Polynomial Differentiated}
We can subject the polynomial to various Cauchy-Euler operator functions. The simplest one is
\begin{equation}
x\partial_{x}{(e^{-x}S_m(x))}=e^{-x}S_{m+1}(x)        \label{eqn13820}
\end{equation}
Repeating this will produce an important relation
\begin{equation}
(x\partial_{x})^j{(e^{-x}S_m(x))}=e^{-x}S_{m+j}(x)        \label{eqn13830}
\end{equation}
We can extend this thinking to an exponential function
\begin{equation}
e^{-\beta{x}\partial_{x}}{(e^{-x}S_m(x))}=\sum^{\infty}_{j=0}{\frac{(-\beta)^j({x}\partial_{x})^j{e^{-x}S_{m}(x)}}{j!}}        \nonumber
\end{equation}
\begin{equation}
=\sum^{\infty}_{j=0}{\frac{(-\beta)^j{e^{-x}S_{m+j}(x)}}{j!}}\equiv{e^{-xe^{-\beta}}S_m(xe^{-\beta})}        \label{eqn13840}
\end{equation}
The last property is achieved by the basic property of the Cauchy-Euler operator in exponential form as presented in Appendix A. We can subject the (RR1) to $e^{-x}$ and $x\partial_x$ to get
\begin{equation}
(x\partial_{x})(e^{-x}S_{m}^{'}(x))=e^{-x}(1-\frac{1}{x})S_{m+1}(x)+\frac{e^{-x}}{x}S_{m+2}(x)          \label{eqn13850}
\end{equation}
\section{Appendix. Integration and Differentiation Results}\label{apx:appc}
\subsection{Operator Inside Integral}
In this section we present some results achieved by integration of various expressions. We may integrate the formal $\zeta$ function operator expression subjected to some integrable function $g(x)$
\begin{equation}
\int_0^{\infty}{\zeta(-\beta{x}\partial_x)g(x)dx}=\sum^{\infty}_{k=1}{\int_0^{\infty}{{k^{\beta{x}\partial_x}}{g(x)dx}}} \label{eqn23300}
\end{equation}
\begin{equation}
=\sum^{\infty}_{k=1}{\int_0^{\infty}{{g(xk^{\beta})dx}}} \label{eqn23310}
\end{equation}
\begin{equation}
=\sum^{\infty}_{k=1}{\frac{\int_0^{\infty}{{g(u)du}}}{k^{\beta}}}=\zeta(\beta)\int_0^{\infty}{{g(u)du}} \label{eqn23320}
\end{equation}
This result can be compared with the M\"{u}ntz formula for the $\zeta(s)$ function
\begin{equation}
\zeta(s)\int_0^{\infty}{y^{s-1}F(y)dy}=\int_0^{\infty}{x^{s-1}\sum^{\infty}_{n=1}{F(nx)dx}} \label{eqn23330}
\end{equation}
\subsection{Operator from the Outside}
Next we have an operator affecting the integral from the outside, getting
\begin{equation}
\zeta(-\beta{x}\partial_x)\int_0^{\infty}{{g(x,y)}h(y)dy}=\sum^{\infty}_{k=1}{\int_0^{\infty}{{h(y)}{g(xk^{\beta},y)dy}}} \label{eqn23340}
\end{equation}
An example follows in the same spirit as above
\begin{equation}
\zeta(-\beta{x}\partial_x)\int_0^{\infty}{{y^{x-1}}F(y)dy}=\sum^{\infty}_{k=1}{\frac{1}{k^{\beta}}\int_0^{\infty}{{F(u^{\frac{1}{k^{\beta}}})}{u^{x-1}du }}} \label{eqn23350}
\end{equation}
\subsection{Integrating the Series}
We can integrate our exponential series $f(x,\alpha)$ with an interesting result
\begin{equation}
\int_0^{x}{{f(x,\alpha)}dx}=\int_0^{x}{{\sum^{\infty}_{k=1}{e^{-tk^{\alpha}}}dt}}=\zeta(\alpha)-\sum^{\infty}_{k=1}{\frac{e^{-xk^{\alpha}}}{k^{\alpha}}} \label{eqn23360}
\end{equation}
We can do it in another way by integrating the factorial series
\begin{equation}
\int_0^{x}{{f(x,\alpha)}dx}=e^{-x}[-\zeta(\alpha)+x\sum^{\infty}_{k=1}{\sum^{\infty}_{m=1}{\frac{(-\alpha{ln(k)})^m}{m!}\sum^{m-1}_{j=1}{(-1)^{j+1}S_j}}}]+\zeta(\alpha) \label{eqn23365}
\end{equation}
which are naturally the same.
\subsection{Integral Operator}
We have the integral operator subjected to the power of $x$
\begin{equation}
(\beta\int{\frac{dx}{x}})^m{x^n}=\frac{\beta^{m}x^n}{n^m}+C \label{eqn23370}
\end{equation}
Integration in the following expression and use of the recursion relation (RR1) will give
\begin{equation}
\int{dx\frac{e^{-x}S_{m+1}(x)}{x}}=e^{-x}S_m(x)+C \label{eqn23400}
\end{equation}
More properties of the polynomial are found in Appendix B.
\subsection{Series Differentiated}
By differentiating the equation (\ref{eqn1080}) we get
\begin{equation}
\partial_x{f(x,\alpha)}=\frac{e^{-x}}{x}\sum^{\infty}_{k=1,m=0}{\frac{(\alpha{ln(k)})^{m}S_{m+1}(x)}{m!}} \label{eqn23590}
\end{equation}
and again
\begin{equation}
\partial_{x}^{2}{f(x,\alpha)}=\frac{e^{-x}}{x^2}\sum^{\infty}_{k=1,m=0}{\frac{(\alpha{ln(k)})^{m}\Delta{S_{m+1}(x)}}{m!}} \label{eqn23600}
\end{equation}
\section{Appendix. Other Series Expressions}\label{apx:appd}
\subsection{The Exponential Function}
The generating function (\ref{eqn2100}) allows solving of a summation formula for the polynomial itself simply by evaluating it at $t=1$
\begin{equation}
e^{-x(e-1)}=\sum^{\infty}_{m=0}{\frac{S_{m}(x)}{m!}}  \label{eqn3890}
\end{equation}
This proves that the infinite sum of the normed polynomial will always converge for a finite $x\in{R}$.

Since we have a new polynomial at hand, we may be able to apply it to get new expressions for common functions. By evaluating the generating function (\ref{eqn2100}) at $t=\pm{i\frac{\pi}{2}}$ or $e^t=\pm{i}$ we have
\begin{equation}
e^{\mp{xi}}=e^{-x}\sum^{\infty}_{m=0}{\frac{(\frac{\pi}{2})^m{S_{m}(x)(\pm{i})^m}}{m!}}  \label{eqn3900}
\end{equation}
or
\begin{equation}
e^{{x(1\mp{}i})}=\sum^{\infty}_{m=0}{\frac{(\frac{\pi}{2})^m{S_{m}(x)(\pm{i})^m}}{m!}}  \label{eqn3901}
\end{equation}
This is a new expression for the exponential function, valid for finite $x$. 
\subsection{The Sine Function}
The result above will give rise to  other series expressions, for instance
\begin{equation}
sin(x)=\frac{e^{-x}}{2i}\sum^{\infty}_{m=0}{\frac{S_{m}(x)(\frac{\pi}{2})^m[(-i)^m-(i)^m]}{m!}}  \nonumber
\end{equation}
\begin{equation}
=-{e^{-x}}\sum^{\infty}_{j=0}{\frac{S_{2j+1}(x)(\frac{\pi}{2})^{2j+1}(-1)^j}{(2j+1)!}}  \label{eqn3910}
\end{equation}
We have illustrated approximation of the function in terms of one and two first terms in Figure \ref{fig:sincosapp} below.
\subsection{The Cosine Function}
\begin{figure}[ht]
	\centering
		\includegraphics[width=1.00\textwidth]{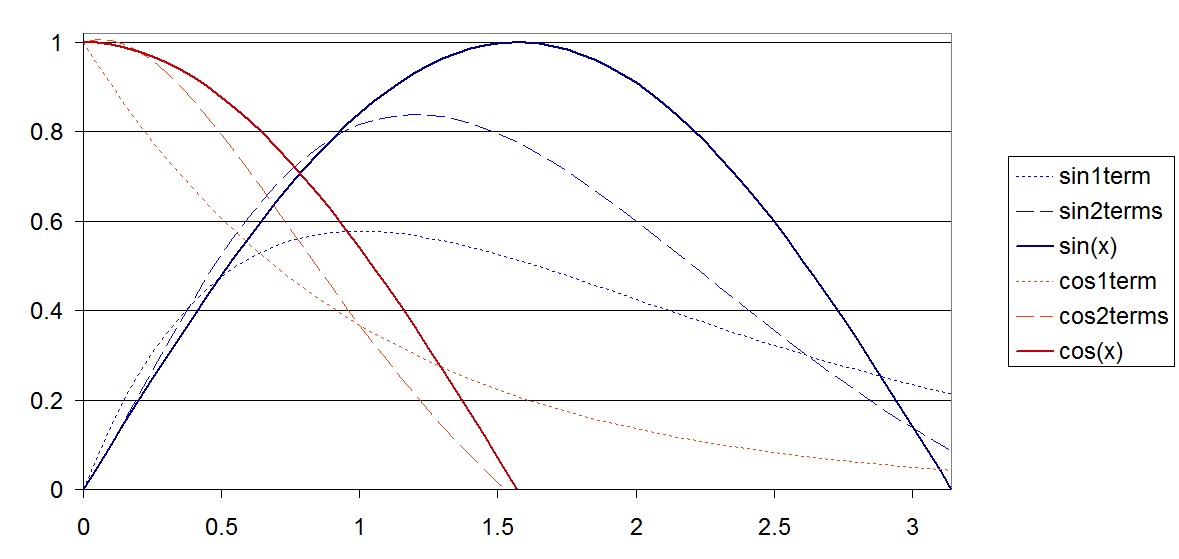}
	\caption{sin(x) and cos(x) approximations with one and two terms compared to accurate functions}
	\label{fig:sincosapp}
\end{figure}

Similarly we can develop equation (\ref{eqn3901}) to
\begin{equation}
cos(x)=\frac{e^{-x}}{2}\sum^{\infty}_{m=0}{\frac{S_{m}(x)(\frac{\pi}{2})^m[(-i)^m+(i)^m]}{m!}}  \nonumber
\end{equation}
\begin{equation}
={e^{-x}}\sum^{\infty}_{j=0}{\frac{S_{2j}(x)(\frac{\pi}{2})^{2j}(-1)^j}{(2j)!}}  \label{eqn3920}
\end{equation}
In these two results above, we have replaced the index with a more suitable one. These series are new representations for the basic trigonometric functions. The range of validity is $x\in$[$0,\frac{\pi}{2})$ which is sufficient for full representation of both functions. The Figure \ref{fig:sincosapp} shows the approximation with one and two first terms of the new series.
\subsection{A Double Sum Identity}
By setting $\alpha=1$ in equation (\ref{eqn2760}) and knowing the result of the simple sum, we obtain a double sum identity
\begin{equation}
\frac{1}{1-e^{-x}}=\sum^{\infty}_{k=1}{\sum^{\infty}_{m=0}{\frac{(ln(k))^{m}S_{m}(x)}{m!}}} 
\end{equation}
By using the generating function relation (GF) we get
\begin{equation}
\frac{1}{1-e^{-x}}=e^{x}\sum^{\infty}_{k=1}{e^{-xe^{ln(k)}}}=e^{x}\sum^{\infty}_{k=1}e^{-xk} 
\end{equation}
which becomes a known identity used before. 
\subsection{Applying the Symmetric Generating Function}
Let's look at the symmetric generating function (\ref{eqn2102}) by setting $y=1$
\begin{equation}
e^{-x}=\sum^{\infty}_{m=0}{\frac{S_{m}(x)(ln(2))^{m}}{m!}} \label{eqn40000}
\end{equation}
or $x=1$ and changing the variable name to $x$
\begin{equation}
e^{-x}=\sum^{\infty}_{m=0}{\frac{S_{m}(1)(ln(x+1))^{m}}{m!}}  \label{eqn40020}
\end{equation}
We can readily apply the relation (\ref{eqn40000}) by assuming it is valid for a complex-valued $z$, at least in the first and fourth quadrants.
\begin{equation}
sin(z)=\frac{1}{2i}\sum^{\infty}_{m=0}{\frac{(S_{m}(-iz)-S_m(iz))}{m!}}
\end{equation}
\begin{equation}
cos(z)=\frac{1}{2}\sum^{\infty}_{m=0}{\frac{(S_{m}(-iz)+S_m(iz))}{m!}}
\end{equation}
These can be repeated yet in another form in terms of the equation (\ref{eqn40020}). These series differ from the usual real-valued series expressions for the trigonometric functions as they consist of both real and imaginary terms due to the polynomials. The imaginary terms cancel each other when all terms are accumulated. While adding terms one by one, the accurate point is approached along a curve which runs through various complex values before finally hitting the accurate real value.

Returning to the symmetric generating function (\ref{eqn2102}) by setting $y=x$ we have the interesting equation for the Gaussian
\begin{equation}
e^{-x^2}=\sum^{\infty}_{m=0}{\frac{S_{m}(x)(ln(x+1))^{m}}{m!}}
\end{equation}
This is shown in graphical form in Figure \ref{fig:exp2} below. The convergence is not too strong. However, when seen on a wider scale (not shown), we can notice that the approximating curves approach closer and closer the correct exponential function when the number of terms is increased. Even the first term alone after the unity
\begin{equation}
e^{-x^2}\approx{1-xln(x+1)}
\end{equation}
gives a fair accuracy while $0<x<1$. Negative values are limited to $-1<x$ but accuracy is not as good on that side. We show the traditional Taylor's expansion approximation with two terms as well. We can see that our new approximation is much better in accuracy in this range.
\begin{figure}[ht]
	\centering
		\includegraphics[width=1.00\textwidth]{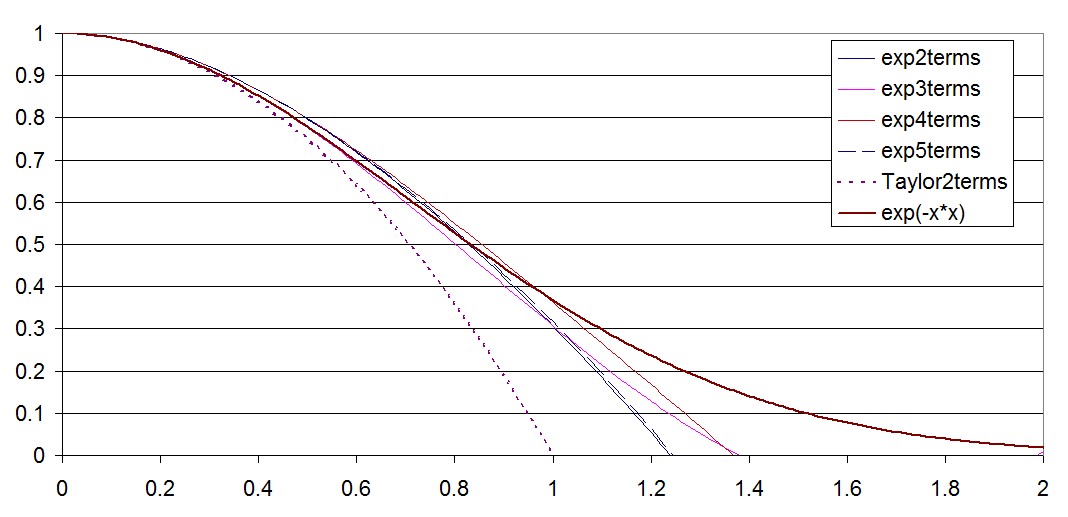}
	\caption{Approximation with a varying number of terms compared to the Gaussian and the Taylor's series with two terms}
	\label{fig:exp2}
\end{figure}
\end{document}